\documentclass[12pt]{article}
\usepackage{amssymb}
\usepackage{amsfonts}
\usepackage{amsmath}
\usepackage[usenames]{color}
\usepackage{mathrsfs}
\usepackage{amsfonts}
\usepackage{amssymb,amsmath}
\usepackage{CJK}
\usepackage{cite}
\usepackage{cases}
\usepackage{amsthm}\usepackage{multirow}
\usepackage{graphicx}
\usepackage{float}
\usepackage{makecell}
\usepackage{latexsym}
\usepackage{CJK}

\def\cl{\centerline}

\def\vs{\vspace*}

\def\ni{\noindent}

\numberwithin{equation}{section}
\newtheorem{theo}{Theorem}[section]

\newtheorem{coro}[theo]{Corollary}
\newtheorem{lemm}[theo]{Lemma}
\newtheorem{exam}[theo]{Example}

\newtheorem{remark}[theo]{Remark}

\begin{document}
\begin{center}
\cl{\large\bf \vs{6pt} Classifications of isoparametric hypersurfaces in}
\cl{\large\bf \vs{6pt} Randers space forms\,{$^*\,$}}
\footnote {$^*\,$ Project supported by NNSFC (No.11471246,11971253), AHNSF (No.1608085MA03) and KLAMFJPU(No. SX201805).
\\\indent\ \ $^\dag\,$ yst419@163.com
}
\cl{Qun He$^1$, Peilong Dong$^1$ and Songting Yin$^{2,3}$$^\dag\,$  }

\cl{\small 1 School of Mathematical Sciences, Tongji University, Shanghai,
200092, China.}
\cl{\small 2 Department of Mathematics and computer science, Tongling University,}
\cl{\small Tongling, 244000, China.}
\cl{\small 3 Key Laboratory of Applied Mathematics (Putian University), } \cl{\small Fujian Province University, Fujian Putian, 351100, China}
\end{center}

{\small
\parskip .005 truein
\baselineskip 3pt \lineskip 3pt

\noindent{{\bf Abstract:}
In this paper, we give the complete classifications of isoparametric hypersurfaces in Randers space forms. By studying the principal curvatures of anisotropic submanifolds in a Randers space~$(N,F)$ with the navigation data~$(h,W)$, we find that a Randers space form $(N, F, d\mu_{BH})$ and the corresponding Riemannian space $(N,h)$ have the same isoparametric hypersurfaces, but in general, their isoparametric functions are different.  We give a necessary and sufficient condition for an isoparametric function of $(N,h)$ to be isoparametric on $(N, F, d\mu_{BH})$, from which we get some examples of isoparametric functions.
\vs{5pt}

\ni{\bf Key words:}
 isoparametric hypersurface, Randers space form, principal curvature, anisotropic submanifold.}

\ni{\it Mathematics Subject Classification (2010):} 53C60, 53C40, 53B25.}
\parskip .001 truein\baselineskip 6pt \lineskip 6pt
\section{Introduction}

In Riemannian geometry, the study of isoparametric hypersurfaces has a long history. Since 1938, E. Cartan had begun to study the isoparametric hypersurfaces in real space forms with constant sectional curvature $c$ systematically.
The classification of isoparametric hypersurfaces in space forms is a classical geometric problem with a history of almost one hundred years. Isoparametric hypersurfaces in Euclidean and hyperbolic spaces were classified in 1930's \cite{C,SF,SS}. For the classification of isoparametric hypersurfaces in a unit sphere, which is the most difficult case, there are many important results (as like \cite{TY, GT}, etc.) and it was recently completely solved in~\cite{C1}.

In Finsler geometry, the concept of isoparametric hypersurfaces has been introduced in~\cite{HYS}. Let~$(N,F,d\mu)$ be an
$n$-dimensional Finsler manifold with volume form~$d\mu$. A function~$f$
on~$(N,F,d\mu)$ is said to be \textit{isoparametric} if  there are~$\tilde a (t)$ and~$\tilde b (t)$ such that
\begin{equation}\label{1.3} \left\{\begin{aligned}
&F(\nabla f)=\tilde{a}(f),\\
&\Delta f=\tilde{b}(f),
\end{aligned}\right.
\end{equation}
where~$\nabla f$ and~$\Delta f$ denote the nonlinear gradient and
Laplacian of~$f$ with respect to $d\mu$, respectively (see Section 2.1 and 2.3 for details).

Studying and classifying isoparametric hypersurfaces in Finsler space forms are interesting problems naturally generalized from Riemannian geometry. In~\cite{HYS}, the authors studied isoparametric hypersurfaces in Finsler space forms, and obtained the Cartan type formula and some classifications on the number of distinct principal curvatures or their multiplicities. For some very special Finsler space forms, such as Minkowski space (with zero flag curvature) and Funk space (a special Rander space with negative constant flag curvature), the isoparametric hypersurfaces have been completely classified~\cite{HYS, GM, HYS1}. Xu \cite{X} studied a special class of isoparametric hypersurfaces in a Randers sphere (with positive constant flag curvature).

Randers metrics play a fundamental role in Finsler geometry. Those with constant flag curvature were classified in~\cite{BRS}, using Zermelo's navigation method. A forward (resp. backward) complete and simply connected Randers manifold with constant flag curvature $c$ is called a \emph{forward (resp. backward) Randers space form}, which is denoted by $(N (c),F)$.
In fact, except in very few cases, the majority of known examples of Finsler space forms with non-zero flag curvature are Randers space forms. So it is natural to consider the isoparametric hypersurfaces in Randers space forms. Unlike the Riemannian case, there are infinitely many Randers space forms, which are not isometric or even are not homothetic to each other. So we cannot expect to solve the problems of classification in a way similar to Riemannian geometry.

In this paper, we will give the complete classifications of isoparametric hypersurfaces in a forward (or backward) Randers space form~$(N^{n}(c),F)$.
 By using navigation process, we find the following
\begin{theo} \label{thm0}
Let~$(N, F)$ be a forward (or backward) Randers space form with the navigation data~$(h, W)$. Then~$(N, F,d\mu_{BH})$ and~$(N,h)$ have the same isoparametric hypersurfaces, and the number of distinct principal curvatures and the multiplicities of each principal curvature are also the same. So the isoparametric hypersurfaces in $(N, F,d\mu_{BH})$ can be completely classified (see Table 1 for the accurate classifications).
\end{theo}

The contents of this paper are organized as follows. In section 2, some fundamental concepts and formulas are given for later use. In section 3,  we consider the principal curvatures of submanifolds with respect to $F$ and $h$, respectively, and derive the classifications of isoparametric hypersurfaces in Randers space forms. Finally, we consider the relation between isoparametric functions with respect to $F$ and $h$ and give some examples of isoparametric functions in special Randers space forms.

\section{Preliminaries}
\subsection{Finsler-Laplacian}
Let~$(N,F)$ be an~$n$-dimensional oriented smooth Finsler manifold and~$TN$ be the tangent bundle over~$N$ with local coordinates
$(x,y)$, where~$x=(x^i)$ and~$y=y^i\frac{\partial}{\partial x^{i}}$. The fundamental form~$g$ of~$(N,F)$ is given by
$$g=g_{ij}(x,y)dx^{i} \otimes dx^{j}, ~~~~~~~g_{ij}(x,y)=\frac{1}{2}[F^{2}] _{y^{i}y^{j}}. ~$$
The projection~$\pi: TN\rightarrow N$ gives rise to the pull-back bundle~$\pi^{\ast}TN$ and its dual bundle
$\pi^{\ast}T^{\ast}N$ over~$TN\backslash{0}$. Recall that on the pull-back bundle~$\pi^{\ast}TN$ there exists a unique
\emph{Chern connection}~$\nabla$ with~$\nabla
\frac{\partial}{\partial x^i}=\omega_{i}^{j}\frac{\partial}{\partial
x^{j}}=\Gamma^{i}_{jk}dx^k\otimes\frac{\partial}{\partial x^{j}}$ satisfying \cite {BCS}
$$dg_{ij}-g_{ik}\omega^{k}_{j}-g_{kj}\omega^{k}_{i}=2FC_{ijk}\delta y^{k},$$
$$\delta y^{i}:=\frac{1}{F}(dy^{i}+N^{i}_{j}dx^{j}),~~~~N_j^i:=\frac{\partial G^i}{\partial y^j}=\Gamma^{i}_{jk}y^k,$$
where $C_{ijk}=\frac{1}{2}\frac{\partial g_{ij}}{\partial y^k}$ is
called the\emph{ Cartan tensor} and \begin{align*}
G^{i}=\frac{1}{4}g^{il}\left\{[F^{2}]_{x^{k}y^{l}}y^{k}-[F^{2}]_{x^{l}}\right\}
\end{align*}
 are the geodesic coefficients of $(N,F)$.

Let~$X=X^{i}\frac{\partial}{\partial x^{i}}$ be a vector field. Then the \emph{covariant derivative} of~$X$ along
$v=v^i\frac{\partial}{\partial x^{i}}\in T_{x}N$ with respect to a reference  vector~$w\in T_{x}N\backslash 0$ for the Chern connection is defined by
\begin{align}
{ D}^{w}_{v}X(x):&=\left\{v^{j}\frac{\partial X^{i}}{\partial x^{j}}(x)+{ \Gamma}^{i}_{jk}(w)v^{j}X^{k}(x)\right\}\frac
{\partial}{\partial x^{i}}.\label{Z1}\end{align}

Let~${\mathcal L}:TN\rightarrow T^{\ast}N$ denote the \emph{Legendre transformation}, satisfying~${\mathcal L}(\lambda
y)=\lambda {\mathcal L}(y)$ for all~$\lambda>0,~y\in TN$, and~\cite{SZ}
\begin{align}\mathcal L(y)&=F(y)[F]_{y^{i}}(y)dx^i,~~\forall y\in TN\setminus \{0\},\\
\mathcal L^{-1}(\xi)&=F^*(\xi)[F^*]_{\xi_{i}}(\xi)\frac{\partial}{\partial x^i},~~\forall \xi\in T^*N\setminus \{0\},\label{Z01}\end{align}
 where~$F^*$ is the dual metric of~$F$. In general,~$\mathcal L^{-1}(-\xi)\neq-\mathcal L^{-1}(\xi)$.
For a smooth function~$f:N\rightarrow \mathbb{R}$, the \emph{gradient vector} of~$f$ at~$x$ is defined as~$\nabla f(x)={\mathcal
L}^{-1}(df(x))\in T_{x}N$.
Set~$N_{f}=\{x\in N|df(x)\neq 0\}$ and~$\nabla^{2}f(x)=D^{\nabla f}(\nabla f)(x)$ for~$x\in N_{f}$. The \emph{Finsler-Laplacian} of~$f$ with respect to the volume form~$d\mu=\sigma(x)dx^{1}\wedge dx^{2}\wedge\cdots\wedge dx^{n}$ is defined by\begin{align}\Delta f=\textmd{div}_{\sigma}(\nabla f)=\textmd{tr}_{g_{\nabla f}}(\nabla^{2}f)-S(\nabla f),
\label{2.22.1}\end{align}
where
\begin{align}\label{1.4.19}
S(x,y)=\frac{\partial G^{i}}{\partial y^{i}}-y^{i}\frac{\partial}{\partial x^{i}}(\ln \sigma(x))
\end{align}
is the \emph{$\mathbf{S}$-curvature}.
\subsection{Anisotropic submanifolds}

Let $(N,F)$ be an~$n$-dimensional Finsler manifold and~$\phi: M\to(N, F)$ be an~$m$-dimensional immersion. Here and from now on, we will use the following convention of index ranges unless other stated:
$$1\leq i, j, \cdots \leq n ;~~~~~~~1\leq a, b, \cdots \leq m< n;$$
$$m+1\leq \alpha, \beta, \cdots  \leq n.$$ Let
$$ \mathcal{N}M=\{(x, \mathbf{n})|~x\in \phi(M),\mathbf{n}\in T_{x}(N), {\mathcal L}(\mathbf{n})(X)=g_{\textbf{n}}(\textbf{n},X)=0,\forall X\in T_xM\},$$
which is called the
{\it normal bundle} of~$\phi$ or~$M$. Note that in general, it is not a vector bundle. We call $\{(M,g_{\textbf{n}})|\textbf{n}\in\mathcal{N}M\}$ an \emph{anisotropic submanifold} of $(N,F)$ to distinguish it from an isometric immersion submanifold $(M,\phi^*F)$.
Moreover, we denote the {\it unit normal bundle} by
$$ \mathcal{N}^0M=\{\mathbf{n}\in \mathcal{N}M|~F(\mathbf{n})=1\}.$$
 For any~$X\in T_xM$ and~$\mathbf{n}\in \Gamma(\mathcal{N}^0M)$, the \emph{shape operator}~${A}_{\mathbf{n}}:T_xM\rightarrow T_xM$ defined by
\begin{equation}\label{0.1}{A}_{\mathbf{n}}(X)=-\left(D^{\mathbf{n}}_{X}\mathbf{n}\right)^{\top}_{g_{\mathbf{n}}}.
\end{equation}
 From Section 3.1 of \cite{HDRY}, we know that ${A}_{\mathbf{n}}$ only
depends on ${\mathbf{n}}(x)$ and
\begin{equation}\label{0.2}g_{\mathbf{n}}({A}_{\mathbf{n}}(X),Y)=g_{\mathbf{n}}(X,{A}_{\mathbf{n}}(Y)).
\end{equation}
The eigenvalues of~${A}_{\mathbf{n}}$,
$\lambda_1,\lambda_2,\cdots,\lambda_{m}$, and~$\hat{H}_{\mathbf{n}}=\sum_{a=1}^{m}\lambda_{a}$ are called the \emph{ principal
curvatures} and the \emph{anisotropic mean curvature} with respect to~$\mathbf{n}$, respectively. If $\lambda_1=\lambda_2=\cdots=\lambda_{m}$ for any~$\mathbf{n}\in \mathcal{N}M$, then $M$ is called to be \emph{anisotropic-totally umbilic}. If~$\hat{H}_{\mathbf{n}}=0$ for any~$\mathbf{n}\in \mathcal{N}M$, then~$M$ is called an \emph{anisotropic-minimal submanifold} of~$(N,F)$.

Let $\phi:M\to N$ be an embedded hypersurface of~$(N,F)$. For any~$x\in M$, there exist exactly two unit\emph{ normal vectors}~$\mathbf{n}_{\pm}$. Let~$\mathbf{n}$ be a given normal vector of
$N$.  Set~$\hat g=\phi^*g_{\mathbf{n}}$.
From \cite{HYS}, we have the following \emph{Gauss-Weingarten formulas}
\begin{align}\label{1.1}{ D}^{\textbf{n}}_{X}Y&={ {\hat{\nabla}}}_{X}Y+\hat{h}(X,Y)\textbf{n},\\
D^{\textbf{n}}_{X}\textbf{n}&=-{A}_{\textbf{n}}X,~~~~~\quad \forall X,~Y\in \Gamma(TM) .\label{1.2}\end{align}
Here
\begin{equation*}\hat{h}(X,Y):=g_{\textbf{n}}(\textbf{n},{ D}^{\textbf{n}}_{X}Y)=\hat {g}({A}_{\textbf{n}}X,Y))\end{equation*}
is called the \emph{second fundamental form}, and
 ${\hat{\nabla}}$ is a torsion-free linear connection on $M$ satisfying \cite{HYS1}
\begin{eqnarray}\label{1.4}
({ \hat{\nabla}}_X\hat g)(Y,Z)=-2C_{\textbf{n}}({A}_{\textbf{n}}X,Y,Z),~~~~~~~\forall X,Y,Z\in \Gamma(TM),
\end{eqnarray}
where $C_{\textbf{n}}$ is the Cartan tensor with $y={\textbf{n}}$.
\subsection{Isoparametric functions and isoparametric hypersurfaces}
 Let $f$ be a non-constant $C^1$ function defined on a Finsler manifold $(N,F,d\mu)$ and smooth on $N_f$. Set $J=f(N_f)$. The function $f$ is said to be \emph{isoparametric} if there exist
a smooth function $\tilde a (t)$ and a continuous function $\tilde b (t)$ defined on $J$ such that (\ref{1.3})
holds on $N_f$.  All the regular level surfaces $M_t = f^{-1}(t)$ form an \emph{isoparametric family}, each of which is called an
\emph{isoparametric hypersurface} in $(N,F,d\mu)$. If $f$ only satisfies the first equation of (\ref{1.3}), then it is said to be {\it transnormal}.

Let $M$ be an embedded hypersurface of $N$. If for any given $x\in M$, there is a neighborhood $V$ of $x$ in $N$ and an isoparametric function $f$ defined on $V$ such that $M\cap V$ is a regular level surface of $f$, then $M$ is called a
\emph{locally isoparametric hypersurface}.
\begin{theo} \label{thm40}
Let $M$ be a  connected and oriented hypersurface embedded in a connected Finsler manifold with constant flag curvature and constant $\mathbf{S}$-curvature. Then $M$ is locally isoparametric if and only if its principal curvatures are all constant.
\end{theo}
\proof  From Theorem 4.2 in \cite{HYS} (which also holds in a local domain), we know that a transnormal function is isoparametric if and only if the principal curvatures of its each regular level surface are all constant.  Then we only need to prove that if the principal curvatures of  $M$ are all constant, then there exists locally a transnormal function $f$ for each point of $M$ such that $M$ is a regular level surface of $f$ and the principal curvatures of each regular level surface of $f$ are all constant.

 Let $\textbf{n}$ be a smooth unit normal vector field of $M$. Because the normal geodesics locally and smoothly depend on $\textbf{n}$, we know that for any $x\in M$, there exist a neighborhood $V$ of $x$ in $N$ and a smooth distance function $f$ defined on $V$ such that $M=f^{-1}(0)$. In fact, we can define $f(p)=d_F(M,p)$ or $f(p)=-d_F(p,M)$ for $p\in V\setminus M$ such that $\nabla f|_M=\textbf{n}$. Then $F(\nabla f)=1$, which shows that $f$ is a transnormal function and its each regular level surface $M_s=f^{-1}(s)$ is a parallel hypersurface of $M$. From Lemma 3.5 in \cite{HDRY}, we know that if the principal curvatures of  $M$ are all constant, then
 the principal curvatures of $M_s$ are also all constant.
\endproof

\section{ Isoparametric hypersurfaces in a Randers space form}
\subsection{Anisotropic submanifolds in Randers spaces}
Let~$(N,F,d\mu_{BH})$ be an~$n$-dimensional Randers space, where
$F=\alpha+\beta=\sqrt{a_{ij}y^iy^j}+b_iy^i$, and let its navigation expression be
$$F=\frac{\sqrt{\lambda h^2+w_{0}^2}-w_{0}}{\lambda}=\frac{\sqrt{\lambda h_{ij}y^iy^j+(w_{i}y^i)^2}-w_{i}y^i}{\lambda},$$
where~$\lambda=1-b^2, ~b=\|W\|_h, ~ W=w^i\frac{\partial}{\partial x^{i}}, ~w_i=h_{ij}w^j$, and
\begin{align}
a_{ij}=\frac{1}{\lambda^2}(\lambda h_{ij}+w_iw_j),~~~~
b_i=-\frac{w_i}{\lambda}.\label{6.00}\end{align} Then
\begin{align}
g_{ij}&=\frac{F}{\alpha}(a_{ij}-\alpha_{y^{i}}\alpha_{y^{j}})+F
_{y^i}F_{y^j}\nonumber\\
&=\frac{F}{\lambda^2\alpha}(\lambda h_{ij}+w_iw_j-\lambda^2\alpha_{y^{i}}\alpha_{y^{j}})+F
_{y^i}F_{y^j}.\label{6.01}\end{align}
Denote the dual metric of $h$ by $h^*$. Then the dual metric of $F$ can be expressed as
\begin{align}
F^*=&h^*+w^0=\sqrt{h^{ij}\xi_i\xi_j}+w^{i}\xi_i,\quad \xi=\xi_idx^i \in T^*_x N.\label{4.1} \end{align}
In fact, the dual metric $F^*$ of a Randers metric $F =\alpha+ \beta$ is still
of Randers type on $T^*M$. If we write $F^*(x; \xi) = \alpha^*(x;\xi) + \beta^*(x; \xi)$, where
$$\alpha^*(x;\xi) =\sqrt{ a^{*ij}(x)\xi_i\xi_j}, \quad \beta^*(x; \xi)= b^{*i}(x)\xi_i, \quad \xi = \xi_idx^i,$$
then
$$a^{*ij} = \frac{(1-\|\beta\|_{\alpha}^2 )a^{ij} + b^ib^j}{(1-\|\beta\|_{\alpha}^2)^2},\quad
b^{*i} = -\frac{b^i }{1-\|\beta\|_{\alpha}^2}.$$
In this case, we can prove that $a^{*ij} = h^{ij}$ and $b^{*i} = w^i$.

Let~$\phi: M\to(N^{n}, F)$ be an $m$-dimensional immersion. Take~$\mathbf{n} \in \mathcal{N}^0(M)$ and~$\nu={\mathcal L}(\mathbf{n})$. From (\ref{Z01}) and (\ref{4.1}), we know that
$$\mathbf{n}^i=F^*_{\xi_i}(\nu)=\frac{h^{ij}\nu_j}{h^*(\nu)}+w^i.$$
Denote~$\bar{\mathbf{n}}=\frac{h^{ij}\nu_j}{h^*(\nu)}\frac{\partial}{\partial x^{i}}$. Then $\bar{\mathbf{n}}$ is a unit normal vector field of~$M$ with respect to~$h$. Thus
\begin{align}
\mathbf{n}=\bar{\mathbf{n}}+W.\label{3.21}
\end{align}
Let $(u^a)=(u^1,\cdots ,u^m)$ be the local coordinates on
$M$ and $d\phi=\phi^i_adu^a\otimes\frac{\partial}{\partial x^{i}}$. Then
\begin{align}
F_{y^{i}}(\mathbf{n})&=\alpha_{y^{i}}(\mathbf{n})-\frac{w_i}{\lambda},\\
\alpha_{y^{i}}(\mathbf{n})\phi^i_a&=F_{y^{i}}(\mathbf{n})\phi^i_a+\frac{w_i\phi^i_a}{\lambda}=\frac{w_i\phi^i_a}{\lambda}.\label{3.220}
\end{align}
It follows from $F(\mathbf{n})=1$ that
\begin{align}\lambda\alpha(\mathbf{n})=\lambda-\lambda\beta(\mathbf{n})=\lambda+\langle {\mathbf{n}}, W \rangle_h=1+\langle\bar{\mathbf{n}}, W \rangle_h.
\label{3.271}
\end{align}
Combing (\ref{6.01}), (\ref{3.21}), (\ref{3.220}) and (\ref{3.271}), yields
\begin{align*}
(\hat g_{\mathbf{n}})_{ab}=&
g_{ij}(\mathbf{n})\phi^i_a\phi^j_b\nonumber\\
=&\frac{1}{\lambda^2\alpha(\mathbf{n})}(\lambda \bar h_{ab}+w_iw_j\phi^i_a\phi^j_b-\lambda^2\alpha_{y^{i}}(\mathbf{n})\alpha_{y^{j}}(\mathbf{n})\phi^i_a\phi^j_b)\nonumber\\
=&\frac{1}{\lambda\alpha(\mathbf{n})} \bar h_{ab},\end{align*}
where $\bar h_{ab}=h_{ij}\phi^i_a\phi^j_b$.
Thus we have the following
\begin{lemm} \label{lem61}Let~$\phi: M\to(N, F)$ be an anisotropic submanifold in a Randers space~$(N,F)$ with the navigation data~$(h,W)$. Then for any smooth section $\mathbf{n}$ of $\mathcal{N}^0(M)$, the induced metric $\hat g_{\mathbf{n}}=\phi^*g_{\mathbf{n}}$ is conformal to $\bar h=\phi^*h$ and satisfies
\begin{align}\hat g_{\mathbf{n}}=\frac{1}{\lambda\alpha(\mathbf{n})} \bar h=\frac{1}{1+\langle \bar{\mathbf{n}}, W \rangle_h}\bar h.
\label{3.270}
\end{align}
\end{lemm}
Denote
$$ r_{ij}=\frac{1}{2}(w_{i|j}+ w_{j|i}),~~~s_{ij}=\frac{1}{2}(w_{i|j}- w_{j|i}),$$
$$ r_{j}=w^{i}r_{ij},~~~r=r_{j}w^{j},~~~s_{j}=w^{i}s_{ij},$$
$$s^{i}=h^{ik}s_{k},~~~r^{i}=h^{ik}r_{k},~~~s^{i}_{~j}=h^{ik}s_{kj},$$
$$s^{i}_{~0}=s^{i}_{~j}y^{j},~~s_{0}=s_{i}y^{i},~~r_{0}=r_{i}y^{i},~~r_{00}=r_{ij}y^{i}y^{j},$$
where~$"|"$ denotes the covariant differential about~$h$.
From \cite{1}, we know that
\begin{align*}
G^{i}=\bar{G}^{i}-Fs^{i}_{~0}-\frac{1}{2}F^{2}(r^{i}+s^{i})+\frac{1}{2}(\frac{y^{i}}{F}-w^{i})(2Fr_{0}-r_{00}-F^{2}r),
\end{align*}
where~$G^{i}$ and~$\bar{G}^{i}$ are the geodesic coefficients of~$F$ and~$h$, respectively. Then
\begin{align}
N^{i}_{j}&=\frac{\partial G^{i}}{\partial y^{j}}=\bar{N}^{i}_{j}-F_{y^{j}}s^{i}_{~0}-Fs^{i}_{~j}-FF_{y^{j}}(r^{i}+s^{i})
\nonumber\\
&+\frac{1}{2}(\frac{\delta^{i}_{j}}{F}-\frac{1}{F^2}y^{i}F_{y^{^{j}}})(2Fr_{0}-r_{00}-F^{2}r)\nonumber\\
&+\frac{1}{2}(\frac{y^{i}}{F}-w^{i})(2F_{y^{j}}r_{0}+2Fr_{j}-2r_{j0}-2FF_{y^{j}}r).\label{3.26}
\end{align}
From \cite{1}, ~$F$ has isotropic~$\mathbf{S}$-curvature if and only if~$W$ satisfies
\begin{align*}
r_{ij}=-2k(x)h_{ij}.
\end{align*}
Using the above formulas, we get
\begin{align}
s^{i}_{~j}=w^{i}_{~|j}+2 k\delta^{i}_{j},~~~r_{j}=-2 k w_{j},\label{3.25}
\end{align}
\begin{align}
r_{0}=-2 k w_{0},~~~r_{00}=-2 k h^{2},~~~r=-2 k b^{2}.\label{3.22}
\end{align}
\begin{lemm} \label{lem6}Let~$\phi: M\to(N, F)$ be an anisotropic submanifold in a Randers space~$(N,F)$ with the navigation data~$(h,W)$. If~$F$ has isotropic~$\mathbf{S}$-curvature~$\mathbf{S}=(n+1)k(x)F$, then for any~$\mathbf{n} \in \mathcal{N}^0(M)$ and~$X\in TM$,
\begin{align}
D^{\mathbf{n}}_{X}\mathbf{n}=\nabla^{h}_{X}\bar{\mathbf{n}}-k(x)d\phi X\label{3.27}.
\end{align}
\end{lemm}
\proof Set $X=X^a\frac{\partial}{\partial u^{a}}$.
By  (\ref{3.21}) $\sim$ (\ref{3.22}), we have
\begin{align*}
D^{\mathbf{n}}_{X}\mathbf{n}=&(\mathbf{n}^{i}_{x^{j}}
+N^{i}_{j}(\mathbf{n}))\phi^j_{a}X^{a}\frac{\partial}{\partial x^{i}}\nonumber\\
=&\left(\mathbf{n}^{i}_{x^{j}}+\bar{N}^{i}_{j}(\mathbf{n})-s^{i}_{~j}+\frac{1}{2}\delta^{i}_{j}(2r_{0}-r_{00}-r)
+\frac{1}{2}(\mathbf{n}^{i}-w^{i})(2r_{j}-2r_{j0})\right)\phi^j_{a}X^{a}\frac{\partial}{\partial x^{i}}\nonumber\\
=&\nabla^{h}_{X}(\bar{\mathbf{n}}+W)-(\nabla^{h}_{X}W+2kd\phi X)+k(-2w_{i}\mathbf{n}^i+h(\mathbf{n})^2+b^2)d\phi X
\nonumber\\
&+k(-2w_{j}\phi^j_{a}X^{a}+2h_{ij}\mathbf{n}^i\phi^j_{a}X^{a}))\bar{\mathbf{n}}\nonumber\\
=&\nabla^{h}_{X}\bar{\mathbf{n}}-2kd\phi X+k(-2\langle \mathbf{n}, W \rangle_h+\langle \mathbf{n}, \mathbf{n} \rangle_h+\|W\|^{2}_{h})d\phi X\nonumber\\
&+2k(-\langle W, d\phi X \rangle_h+\langle \mathbf{n}, d\phi X\rangle_h)\bar{\mathbf{n}}\nonumber\\
=&\nabla^{h}_{X}\bar{\mathbf{n}}-2k d\phi X+k |\bar{\mathbf{n}}|^{2}d\phi X\nonumber\\
=&\nabla^{h}_{X}\bar{\mathbf{n}}-k d\phi X.
\end{align*}
\endproof
Thus, we have the following
\begin{theo} \label{thm6}
Let~$M$ be an anisotropic submanifold in a Randers space~$(N,F,d\mu_{BH})$ with the navigation data~$(h,W)$. If~$F$ has isotropic~$\mathbf{S}$-curvature~$\mathbf{S}=(n+1)k(x)F$, then for any $\mathbf{n} \in \mathcal{N}^0(M)$, the shape operators of $M$ in Randers space~$(N, F)$ and Riemannian space~$(N, h)$, ${A}_{\mathbf{n}}$ and $\bar{A}_{\bar{\mathbf{n}}}$,  have the same principal vectors and their principal curvatures satisfy
\begin{align}\label{3.24}\lambda=\bar{\lambda}+k(x),\end{align}
where~$\lambda$ and~$\bar{\lambda}$ are the principal curvatures of~$M$ in Randers space~$(N, F)$ and Riemannian space~$(N, h)$, respectively.\end{theo}
\proof Set $X=X^a\frac{\partial}{\partial u^{a}}$ and $\phi_a=d\phi\frac{\partial}{\partial u^{a}}$.
By (\ref{0.1}) and (\ref{3.27}), we know that
\begin{align*}
{A}_{\mathbf{n}}X=-\left[D_{X}^{\mathbf{n}}\mathbf{n}\right]^{\top}_{g_{\mathbf{n}}}
=-g_{\mathbf{n}}(\nabla^{h}_{X}\bar{\mathbf{n}},\phi_a)(\hat g_{\mathbf{n}})^{ab}\frac{\partial}{\partial u^{b}}+kX.
\end{align*}
From (\ref{6.00})$\sim$(\ref{3.27}) and (\ref{1.4}), we have

\begin{align*}
-g_{\mathbf{n}}(\nabla^{h}_{X}\bar{\mathbf{n}},\phi_a)(\hat g_{\mathbf{n}})^{ab}\frac{\partial}{\partial u^{b}}
=&-\frac{1}{\lambda^2\alpha(\mathbf{n})}(\lambda \bar h_{ij}+w_iw_j-\lambda^2\alpha_{y^{i}}(\mathbf{n})\alpha_{y^{j}}(\mathbf{n}))\phi^i_a\bar{\mathbf{n}}^j_{|c}X^c(\hat g_{\mathbf{n}})^{ab}\frac{\partial}{\partial u^{b}}\nonumber\\
=&-\frac{1}{\lambda}\left(\lambda \bar h_{ij}\phi^i_a+w_i\phi^i_a(w_j-\lambda\alpha_{y^{j}}(\mathbf{n})\right)\bar{\mathbf{n}}^j_{|c}X^c\bar h^{ab}\frac{\partial}{\partial u^{b}} \nonumber\\
=&-\left[\nabla^{h}_{X}\bar{\mathbf{n}}\right]^{\top}_{h}+\langle \phi_a, W \rangle_hF_{y^{j}}(\mathbf{n})\bar{\mathbf{n}}^j_{|c}X^c\bar h^{ab}\frac{\partial}{\partial u^{b}} \nonumber\\
=&\bar{A}_{\bar{\mathbf{n}}}X+(W)^{\top} _hg_{\mathbf{n}}(\mathbf{n},\nabla^{h}_{X}\bar{\mathbf{n}}) \nonumber\\
=&\bar{A}_{\bar{\mathbf{n}}}X+(W)^{\top} _hg_{\mathbf{n}}(\mathbf{n},D^{\mathbf{n}}_{X}{\mathbf{n}}) \nonumber
=\bar{A}_{\bar{\mathbf{n}}}X .
\end{align*}
Thus ${A}_{\mathbf{n}}X=\bar{A}_{\bar{\mathbf{n}}}X+k(x) X$, which shows ${A}_{\mathbf{n}}$ and $\bar{A}_{\bar{\mathbf{n}}}$ have the same principal vectors and the proof is completed.
\endproof

The following corollaries are immediate consequences of Theorem 3.3, so we skip their proofs.
\begin{coro}In a Randers space~$(N,F,d\mu_{BH})$ with isotropic~$\mathbf{S}$-curvature and the navigation data~$(h,W)$, $M$ is anisotropic-totally umbilic if and only if it is totally umbilic in Riemannian space~$(N, h)$.
\end{coro}
\begin{coro}In a Randers space~$(N,F,d\mu_{BH})$ with constant~$\mathbf{S}$-curvature and the navigation data~$(h,W)$, the principal curvatures of an anisotropic submanifold are all constant if and only if its principal curvatures in Riemannian space~$(N, h)$ are all constant.
\end{coro}
\begin{coro}In a Randers space~$(N,F,d\mu_{BH})$ with constant~$\mathbf{S}$-curvature and the navigation data~$(h,W)$, an anisotropic submanifold $M$ has constant mean curvature if and only
if $M$ also has constant mean curvature in Riemannian space $(N, h)$. Especially, when the $\mathbf{S}$-curvature vanishes, $M$ is anisotropic minimal if and only if it is minimal in Riemannian space~$(N, h)$.
\end{coro}
\subsection{Classifications of isoparametric hypersurfaces in a Randers space form}
\textbf{{Proof of Theorem \ref{thm0}}:}
\proof Let~$(N, F,d\mu_{BH})$ be a Randers space with the navigation data~$(h, W)$. By Theorem 5.11 in \cite{1}, $(N, F)$ has constant flag curvature~$c$ if and only if the Riemannian space~$(N, h)$ has constant sectional curvature~$\bar c $ and~$W$ is a homothetic vector field with dilation~$k_0$. In this case,~$F$ has constant~$\mathbf{S}$-curvature, that is, ~$k(x)=k_0~$ in (\ref{3.24}) and $c=\bar c -k_0^2$. Then Theorem \ref{thm0} follows from Theorem \ref{thm40} and Theorem \ref{thm6}.
\endproof
  Note that in general, when $(N(c), F)$ is complete, $(N(\bar c), h)$ is not necessarily complete. In this case, $(N(\bar c), h)$ is isometric to an open subset of a real space form $\bar N(\bar c)$. That is $$N(\bar c)\cong \{x\in \bar N(\bar c)~|~\|W\|_h<1\},$$ where $\bar N(\bar c)=\mathbb{R}^{n},~\mathbb{H}^{n}(\bar c)$ or $\mathbb{S}^{n}(\frac{1}{\sqrt{\bar c}}).$
By Proposition 5.4 in \cite{1}, $k_0=0$ when $\bar c=c+k_0^2\neq0$. Specifically, we have
\\
 $(1)$ if $c=0$, then $\bar c=k_0=0$, and $N(\bar c)\cong \{x\in \mathbb{R}^{n}~|~\|W\|_h<1\}$;\\
 $(2)$ if $c<0$, then $\bar c=c+k_0^2=0$ for $k_0\neq0$, and $N(\bar c)\cong \{x\in \mathbb{R}^{n}~|~\|W\|_h<1\}$; or $\bar c=c<0$ for $k_0=0$, and $$N(\bar c)\cong \{x\in \mathbb{H}^{n}(c)~|~\|W\|_h<1\};$$ \\ $(3)$ if $c>0$, then $(N(c), F)$ is compact and $\bar c=c$, which implies that $(N(\bar c),h)\cong \mathbb{S}^{n}(\frac{1}{\sqrt{c}}).$ Thus the Randers space form $(N(c), F)$ is globally isometric to a Randers sphere.

From the works of Cartan \cite{C, C40} and M\"{u}nzner \cite{M80}, we know that any connected
isoparametric hypersurface embedded in a real space form is contained in a unique
complete isoparametric hypersurface. Up to now, the classifications of the isoparametric hypersurfaces in real space forms have been completely solved \cite{C, SF, SS, C1, CS1, C2, C39, C40,M80, OT76, C3, DN85, M13, FKM81}. So according to Theorem 1.1, we can give the complete classifications of isoparametric hypersurfaces in a Randers space form~$(N(c), F)$. The classification results are summarized in the following  table.
\begin{table}[h]
  \begin{center}
  \caption{Classification results for isoparametric hypersurfaces in $(N(c), F, d\mu_{BH})$}
    \begin{tabular}{|c|c|c|c|c|c|l|c|}
  \hline
     $K_F=c$ &$\mathbf{S}$-curv. & $N(c)$ &$g$ &dim$M$ & mul.& $M$ is an open subset &main \\
      &&&&&&of following hypersurfaces &ref. \\
  \cline{1-8}
      $c=0$ &$k_0=0$  & $\mathbb{R}^{n}$ &$g$=1& $n$-1&$n$-1 & a hypersphere $\mathbb{S}^{n-1}$&\multirow{4}{*}{}\\
      &&$\|W\|_h<1$&&&& or a hyperplane $\mathbb{R}^{n-1}$&\cite{C}\\
 \cline{1-2} \cline{4-7}
     \multirow{2}{*}{$c=-1$}  &$k_0^{2}=1$ &&$g$=2& $n$-1&($m$,$n$-$m$-1) &a cylinder $\mathbb{S}^{m}\times$ $\mathbb{R}^{n-m-1}$ &\cite{SF}\\
  \cline{2-7}
       & \multirow{2}{*}{$k_0=0$}&$\mathbb{H}^{n}$&$g$=1& $n$-1&$n$-1 & a sphere $\mathbb{S}^{n-1}$, a hyperbolic &\cite{SS}\\
       &&$\|W\|_h<1$&&&& $\mathbb{H}^{n-1}$ or a horosphere $\mathbb{R}^{n-1}$&\\
  \cline{4-7}
      &&& $g$=2&$n$-1&($m$,$n$-$m$-1)& a cylinder $\mathbb{S}^{m}\times \mathbb{H}^{n-m-1}$ &\\
  \cline{1-8}
       \multirow{6}{*}{\emph{$c=1$}}& \multirow{6}{*}{\emph{$k_0=0$}}&\multirow{6}{*}{$\mathbb{S}^{n}$}& $g$=1 &$n$-1 &$n$-1& a great or small &\\
       &&&&&& hypersphere &\\
  \cline{4-7}
      &&& $g$=2&$n$-1&($m$,$n$-$m$-1)&a Clifford torus &\cite{C40}\\
      &&&&&&$S^{m}(r)\times S^{n-m-1}(s)$, &\cite{C39}\\
      &&&&&& $r^{2}+s^{2}=1$&\\
  \cline{4-7}
      &&& \multirow{4}{*}{$g$=3}&3&(1,1)&a tube over a standard&\\
  \cline{5-6}
      &&$\|W\|_h<1$&&6&(2,2)& Veronese embedding  of&\\
  \cline{5-6}
      &&&&12&(4,4)&  $\mathbb{F}$P into $S^{3m+1}$,  where &\\
      &&&&&&$\mathbb{F}=\mathbb{R}$,$\mathbb{C},\mathbb{H}$ or $\mathbb{O}$, for &\\
  \cline{5-6}
      &&&&24&(8,8)& $m=1, 2, 4, 8$, respectively. &\\
 \cline{4-8}
      &&&\multirow{4}{*}{$g$=4}&$2(m_{1}+m_{2})$&$(m_{1},m_{2})$&\multirow{3}{*}{OT-FKM type or}  &\cite{CS1}\\
      &&&&$m_{2}\geq2m_{1}-1$&&&\cite{FKM81}\\
  \cline{5-6}
      &&&&8&(2,2)&homogeneous &\cite{M80}\\
  \cline{5-6}
      &&&&18&(4,5)&&\cite{OT76}\\
  \cline{5-6}
      &&&&14&(3,4)&&\cite{C2}\\
      \cline{5-6}
      &&&&30&(6,9)&&\cite{C3}\\
  \cline{5-6}
      &&&&30&(7,8)&&\cite{C1}\\
  \cline{4-8}
   &&&\multirow{2}{*}{$g$=6}&6&(1,1)&\multirow{2}{*}{homogeneous} &\cite{DN85}\\
  \cline{5-6}
      &&&&12&(2,2)& &\cite{M13}\\
      \hline
    \end{tabular}
  \end{center}
\end{table}
\begin{remark}
For $g=4$, there exist homogeneous examples with multiplicities as follows ($m\geq2$): (1, m-1), (2, 2m-1), (4, 4m-5), (2, 2), (4, 5), (6, 9). Especially, (2, 2) and (4, 5) are just the homogeneous cases that are not included in OT-FKM type.
\end{remark}
\begin{remark}
The isoparametric hypersurfaces constructed by Ferus, Karcher and M$\ddot{u}$nzner are usually referred to as isoparametric hypersurfaces of FKM type \cite{FKM81}. These examples are a generalization of the work of Ozeki and Takeuchi \cite{OT76}, so they are sometimes collectively called OT-FKM type.
\end{remark}
\subsection{ Isoparametric functions of  Randers space forms}
 We know from Theorem \ref{thm0} that the Randers space form~$(N, F, d\mu_{BH})$ with the navigation data~$(h, W)$ and the Riemannian space form~$(N, h)$ have the same isoparametric hypersurfaces, but in general, they have different isoparametric families. For example, from Table 1, it is easy to see that every Euclidean sphere is isoparametric in Minkowski-Randers space or Funk space~$(N^n,F)$.
  But in general, its isoparametric family in~$(N, F)$ is a family of spheres with variational center \cite{HYS, HYS1},  while its isoparametric family in~$(N, h)$ is a family of concentric spheres. Thus the corresponding isoparametric functions are different. On the other hand, the metrics on an isoparametric hypersurface induced from Randers metric~$F$ and Riemannian metric~$h$ are different, and their geometric characteristics are different accordingly. It is still necessary to study the properties of isoparametric functions in Randers space forms.
  \begin{theo} \label{thm02}
 Let~$(N, F,d\mu_{BH})$ be a Randers space with the navigation data~$(h,W)$, where ~$W$ is a homothetic vector field, and let~$f$ be an isoparametric function of~$(N,h)$. Then~$f$ is an isoparametric function of~$(N, F)$ if and only if there is a smooth function~$\varphi$ such that~$df(W)=\varphi(f)$.
\end{theo}
\proof From \cite{HYS1}, we know that  in a Randers space $(N, F, d\mu)$ with the navigation data~$(h, W)$, $f$ is an isoparametric function if and only if~there exist two functions~$\tilde a (t)$ and~$\tilde b (t)$ such that $f$ satisfies
\begin{equation}\label{3.10}\left\{\begin{aligned}
      &|df|_h+\langle df,W^*\rangle_h=\tilde a(f),\\
      &\frac{1}{ |df|_h}\Delta^h_{\sigma}f+\text{div}_{\sigma} W +\frac{1}{ |df|_h^2}\langle d\langle df,W^*\rangle_h,df\rangle_h=\frac{\tilde b(f)}{\tilde a(f)}.
\end{aligned}\right.\end{equation}
 If $f$ is an isoparametric function of~$(N, h)$, then there exist two functions $ a (t)$ and~$ b (t)$ such that $f$ satisfies \begin{equation}\label{3.11}\left\{\begin{aligned}
      &|df|_h=a(f),\\
      &\Delta^h_{\sigma}f=b(f).
\end{aligned}\right.\end{equation}
Furthermore, if $f$ is also an isoparametric function of~$(N, F)$, then the desired conclusion follows directly from the first equation of (\ref{3.10}).

 Conversely, suppose there is a smooth function~$\varphi: f(M)\rightarrow R$ such that~$df(W)=\varphi(f)$. Since $W$ is a homothetic vector
field, we have
$$\text{div}_{\sigma} W=h^{ij}w_{i|j}=-2nk_0,$$
$$\langle d\langle df,W^*\rangle_h,df\rangle_h=\varphi'(f)|df|_h^2.$$ From the above formulas and (\ref{3.11}),  we obtain that
\begin{equation*}\left\{\begin{aligned}
      &|df|_h+\langle df,W^*\rangle_h=a(f)+\varphi(f),\\
      &\frac{1}{ |df|_h}\Delta^h_{\sigma}f+\text{div}_{\sigma} W +\frac{1}{ |df|_h^2}\langle d\langle df,W^*\rangle_h,df\rangle_h=\frac{b(f)}{ a(f)}-2nk_0+\varphi'(f).
\end{aligned}\right.\end{equation*}
Then by (\ref{3.10}),~$f$ is an isoparametric function of~$(N,F)$ with
\begin{equation*}\left\{\begin{aligned}
      &\tilde a(f)=a(f)+\varphi(f),\\
      &\tilde b(f)=b(f)-a(f)\left(2nk_0-\varphi'(f)\right).
\end{aligned}\right.\end{equation*}
\endproof
\begin{remark}
In~\cite{X}, Xu obtained this result for a special case, $df(W)=0$, in a different way.
\end{remark}
By Theorem \ref{thm02}, we can find some examples of isoparametric functions in special Randers space forms.

  Let
$\widetilde{F}$ be a Randers metric with the navigation data~$(\tilde{h}, \widetilde{W})$, where $\tilde{h}=\sqrt{\sum_{\alpha}(y^\alpha)^2},~\widetilde{W}=xQ+x_0,~x,y\in \mathbb R^{n+1},$  $Q$ is an antisymmetric matrix, and $x_0 \in \mathbb R^{n+1}$ is a constant vector. Set $\widetilde{M}=\{x\in \mathbb R^{n+1}~|~|xQ+x_0|^2<1\}$.  Then $\widetilde{W}$ is a Killing vector field in $\mathbb R^{n+1}$ and thus $(\widetilde{M},\tilde{F})$  has constant flag curvature $c=0$. But it is not a local Minkowski space.
\begin{exam}\label{eg.1} Take$$Q=\left(
                      \begin{array}{cc}
                        Q' & 0 \\
                        0 & 0 \\
                      \end{array}
                    \right).
$$Then $\Phi(x)=\langle x,e_{n+1}\rangle$ is an isoparametric function of $(\widetilde{M},\tilde{h})$ with $g=1$, where~$e_{n+1}=(0,\ldots,0,1)$, and
$$\langle\nabla^{\tilde{h}}\Phi,W\rangle=\langle e_{n+1},xQ+x_0\rangle=\langle e_{n+1},x_0\rangle.$$
From Theorem \ref{thm0} and Theorem \ref{thm02}, we know that~$\Phi$ is an isoparametric function of Randers space~$(\widetilde{M},\tilde{F})$ with $g=1$.
\end{exam}
Let~$(\mathbb S^n,h)\hookrightarrow \mathbb R^{n+1} (n\geq2)$ be the standard Euclidean sphere. Take $\widetilde{W}=x Q$, then~it is easy to prove that $W=\widetilde{W}|_{\mathbb S^{n}}$ is a Killing vector field on $\mathbb S^{n}$. Let $F$  be a Randers metric on $\mathbb S^n$ with the navigation data~$(h,W)$. Then $(\mathbb S^{n},{F})$  has constant flag curvature $c=1$, which means  $(\mathbb S^{n},{F})$ is a Randers sphere.
\begin{exam} In Example \ref{eg.1},
take $x_0=0$ and $f=\Phi|_{\mathbb S^{n}}$. Then~$f$ is an isoparametric function of $(\mathbb S^n,h)$ with $g=1$ and
~$$\langle\nabla^hf,W\rangle_h=\langle\nabla^E\Phi-\Phi x,W\rangle=\langle e_{n+1}-\Phi x,xQ\rangle=0.$$
From Theorem \ref{thm0} and Theorem \ref{thm02}, we know that~$f$ is also an isoparametric function of Randers sphere~$(\mathbb S^n,F)$ with $g=1$.\\
\end{exam}
\begin{exam}  Take~$$Q=\left(
                      \begin{array}{cc}
                        Q_1 & 0 \\
                        0 & Q_2 \\
                      \end{array}
                    \right).
$$ Define ~$$\Phi:\mathbb R^{m}\times\mathbb R^{n-m+1}\rightarrow \mathbb R~~~~~~~~$$
$$~~~~~~~(x_1,x_2) ~~~\mapsto ~~~|x_1|^2-|x_2|^2.$$
 Then from \cite{CR}, $f=\Phi|_{\mathbb S^{n}}$ is an isoparametric function of $(\mathbb S^n,h)$ with $g=2$ and
$$\langle\nabla^hf,W\rangle_h=\langle\nabla^{\tilde{h}}\Phi-2\Phi x,W\rangle=\langle 2(x_1,-x_2)-2\Phi x,(x_1,x_2)Q\rangle=0.$$
By Theorem \ref{thm0} and Theorem \ref{thm02}, we conclude that~$f$ is also an isoparametric function of Randers sphere~$(\mathbb S^n,F)$ with $g=2$.
\end{exam}

Qun He \\
School of Mathematical Sciences, Tongji University, Shanghai, 200092, China\\
E-mail: hequn@tongji.edu.cn \\

Peilong Dong\\
School of Mathematical Sciences, Tongji University, Shanghai, 200092, China\\
E-mail: 1710384@tongji.edu.cn\\

Songting Yin\\
Department of Mathematics and Computer Science, Tongling University, Tongling 244000,
China\\
E-mail: yst419@163.com
\end{document}